\documentclass[times,10pt]{article}

\begin{document}

\newtheorem{Remark}{Remark}[section]
\newtheorem{Lemma}{Lemma}[section]
\newtheorem{Proposition}{Proposition}[section]
\newtheorem{Theorem}{Theorem}[section]
\newtheorem{Exercise}{Exercise}[section] 
\newtheorem{Corollary}{Corollary}[section] 
\newtheorem{Definition}{Definition}[section]
\newtheorem{Example}{Example}[section]
\def\ad{\mbox{ad\,}} \def\tr{\mbox{tr\,}}\def\char{\mbox{char\,}}
\def\mod{\mbox{mod\,}}\def\Ad{\mbox{Ad\,}}
\def\diag{\mbox{diag\,}}\def\ad{\mbox{ad\,}} \def\tr{\mbox{tr\,}} 
\def\End{\mbox{End\,}}\def\GL{\mbox{GL\,}}\def\SL{\mbox{SL\,}}
\def\Out{\mbox{Out\,}}\def\Int{\mbox{Int\,}}
\def\Mat{\mbox{Mat\,}}\def\Hom{\mbox{Hom\,}}\def\Iso{\mbox{Iso\,}}
\def\Aut{\mbox{Aut\,}}\def\gl{\mbox{gl\,}}\def\sl{\mbox{sl\,}}
\def\o{\mbox{o\,}}\def\sp{\mbox{sp\,}}
\def\im{\mbox{im\,}}\def\ker{\mbox{ker\,}}\def\deg{\mbox{deg\,}}
\def\id{\mbox{id\,}}\def\mod{\mbox{mod\,}}
\def\grad{\mbox{\rm grad\,}}\def\rot{\mbox{\rm rot\,}}
\def\div{\mbox{\rm div\,}}\def\Grad{\mbox{\rm Grad\,}}
\def\det{\mbox{\rm det\,}} \def\ctg{\mbox{\rm ctg\,}}\def\tg{\mbox{\rm tg\,}}
\def\sn{\mbox{\rm sn\,}}\def\th{\mbox{\rm th\,}}\def\dn{\mbox{\rm dn\,}}
\def\pd#1,#2{\frac{\partial#1}{\partial#2}}\def\diag{\mbox{\rm diag\,}}
\def\v#1{\overline{#1}}\def\sh{\mbox{\rm sh\,}}\def\ch{\mbox{\rm ch\,}}
\def\d#1,#2{\frac{d#1}{d#2}}\def\th{\mbox{\rm th\,}}
\def\qed{\hbox{${\vcenter{\vbox{
\hrule height 0.4pt\hbox{\vrule width 0.4pt height 6pt
\kern 5pt\vrule width 0.4pt}\hrule height 0.4pt}}}$}}

\begin{titlepage} 
\title{Universal Lie algebra extensions via commutative structures} 
\author{A  B  Yanovski \thanks{\it On leave of absence from the Faculty of
Mathematics and Informatics, St.  Kliment Ohridski University, James
Boucher Blvd, 1164 - Sofia, Bulgaria}}

\maketitle 
\begin{center} 
Universidade Federal de Sergipe, Cidade Universit\'aria,\\ ''Jos\'e
Alo\'{\i}sio de Campos'', 49.100-000-S\~ao Christ\'ov\~ao, SE-Brazil
\end{center} 
\begin{abstract} 
We consider some special type extensions of an arbitrary Lie algebra,
which we call universal extensions. We show that these extensions are
in one-to-one correspondence with finite dimensional associative
commutative algebras. We also construct a special kind of these
extensions, that correspond to a finite commutative monoids.
\end{abstract}
\begin{center} 
37K30,  37K05,  17B80
\end{center} 

\end{titlepage}
\section{Introduction}
Let ${\cal G}$ be an arbitrary Lie algebra over the field ${\bf K}$.
Usually, having in mind the applications, ${\bf K}$ is considered to
be one of the classical fields ${\bf R}$ or ${\bf C}$ but there is no
immediate necessity to limit ourself. Let ${\cal G}$ be a Lie
algebra over ${\bf K}$ and let us consider the linear space ${\cal
G}^n$ with elements
\begin{equation}
{\bf x}=(x_1,x_2,\ldots,x_n),\quad x_i\in {\cal G}
\end{equation}
In \cite{JLT}, there was propounded the idea to study the 
possibilities of defining Lie brackets over ${\cal G}^n$, using the
brackets in ${\cal G}$, in such a way that the construction will be
the same for all ${\cal G}$. For this reason we call these brackets
universal ones and the corresponding Lie algebra structures universal
extensions. More specifically, let us introduce brackets obeying the
following properties:
\begin{itemize}
\item The brackets have the form:
\begin{equation}\label{eq:DEF}
([{\bf x},{\bf y}]_W)_s=\sum\limits_{i,j=1}^n W^{ij}_s[x_i,y_j];
\quad W^{ij}_s\in {\bf K}
\end{equation}
\item For $W=(W^{ij}_s)$ fixed, (\ref{eq:DEF}) is a Lie bracket for
arbitrary Lie algebra ${\cal G}$. 
\end{itemize}
For obvious reasons this algebra will be denoted by ${\cal G}^n_W$.
Let us also remark that $W=(W^{ij}_s)$ has a properties of a
$(2,1)$ tensor over ${\bf K}^n$.

The existence and classification of the above structures in the case
${\bf K}={\bf C}$ was considered in
\cite{JLT}, some additional examples of such structures are given in
\cite{ABY}. However the viewpoint in \cite{ABY} is completely
different, we tried there to establish a connection with the
so-called bundles of Lie algebras, see \cite{Ya3,TrFo}. As to the
interest in the universal extensions, it is motivated by the study of a
number of Hamiltonian structures for dynamic systems related
with Mechanics, Hydrodynamics, Magnetohydrodynamics as well as in
other areas, see
\cite{JLT} for an extensive bibliography, and the possibility to
interpret the corresponding Poisson tensors for them as Kirillov
tensors on ${\cal G}_W^*$. Recently, see \cite{JLT1}, an interesting mechanical
system having the Kirillov structure resulting from one of the most
simple extensions has been considered. (For the definitions and the
applications of the Kirillov tensors see \cite{Ki,Lich,MR}). The
structure of the universal extensions is also useful in the questions
about the stability, as there is a possibility to calculate the
Casimir invariants, see again \cite{JLT,JLT1}.

The bilinear form $[{\bf x},{\bf y}]_W$ is a Lie
bracket of the type we are looking for if the tensor $W=(W^{ij}_s)$
satisfies the conditions:
\begin{equation}\label{eq:0}
\begin{array}{c}
W^{ij}_s=W^{ji}_s\\[4pt]
\sum\limits_{k=1}^n(W_i^{sk}W_{k}^{qp}-W_i^{qk}W_{k}^{sp})=0
\end{array}
\end{equation}
In the case ${\bf K}={\bf C}$ the tensors $W_{s}^{ij}$ and the
corresponding Lie algebra extensions have been analyzed in
\cite{JLT}. It was be shown that they are divided into two classes,
called semisimple extensions and solvable extensions.  Knowing the
extensions of the second class we generate the first one and
vice-versa, so it is natural to study only the solvable extensions.
Let us briefly outline the main points of the mentioned analysis and
introduce the requisite definitions.

Instead of studying the tensors $W^{ij}_k$ it is more convenient to
study an equivalent object - the set of matrices: $W^{(i)}$, $i=1,2,\ldots, n$,
with components $(W^{(i)})^j_k=W^{ij}_k$. There is indeed a strong
reason for this, because the second equation in (\ref{eq:0})
actually means that the matrices $W^{(i)}$ commute. As a consequence
in the case ${\bf K}={\bf C}$ by a similarity transformation
$X\mapsto A^{-1}XA$ defined by a nonsingular matrix $A=(A_i^{i'})$
all the ${W}^{(s)}$ can be put simultaneously into a block-diagonal
form, each block being low-triangular with the corresponding generalized
eigenvalue on the diagonal. The linear transformation $A$ applied to
the tensor $W_k^{ij}$ of course gives
\begin{equation}\label{eq:trm1}
W^{i'j'}_{k'}=\sum\limits_{i,j,k=1}^n A^{i'}_iA^{j'}_jA^k_{k'}W_{k}^{ij}
\end{equation}
where as usual $\sum\limits_{k'}A_k^{k'}A_{k'}^s=\delta_k^s$. As a
result the set $\{W^{(i)}\}_{i=1}^n$ transforms as follows:
\begin{equation}\label{eq:trm2}
W^{(i')}=\sum\limits_{i=1}^n(A^{-1}W^{(i)}A)A^{i'}_i 
\end{equation}
and the block structure already obtained is preserved. The block
structure corresponds to a splitting of the algebra ${\cal 
G}_W^{n}$ into a direct sum. Therefore we can limit ourselves with the
irreducible case and can assume that all $W^{(s)}$ are
low triangular.  Now, the symmetry of $W_{i}^{jk}$ entails that the
generalized eigenvalues of $W^{(s)}$ for $s>1$ are zero and so for
$s>1$ the matrices $W^{(s)}$ are low-triangular with zeroes on the
diagonal, and hence are nilpotent.  All diagonal elements of
$W^{(1)}$ can be assumed to be equal to one and the same number
$a\neq 0$  or all to be equal to $0$. The first case is
called semisimple case and the second - solvable. In what follows we
shall ordinary consider solvable extensions and so for all $i$
$(W^{(i)})^j_k=0$ for $j\geq k$. In particular, $W^{(n)}$ is always
the zero matrix. In the semisimple case, as shown in \cite{JLT}, one
can make a linear transform after which $W^{(1)}={\bf 1}$, conserving
the block-triangular form of the matrices $W^{(i)}$ and the fact that
the diagonal elements of $W^{(j)}$ for $j>1$ are zero.

The lower triangular form of the matrices $W^{(i)}$
following from the above discussion in the case ${\bf K}={\bf C}$
will be called canonical form. In terms of the matrices $W^{(i)}$ the
process of passing from semisimple to solvable extension looks in the
following way. In the set $\{W^{(i)}\}_{i=1}^n$ of $n\times n$
matrices we just drop $W^{(1)}$ and form a set
$\{R^{(s)}\}_{s=1}^{n-1}$ where $R^{(s)}$ is constructed from
$W^{(s+1)}$ cutting off the first row and column.  For this reason,
in the semisimple case is useful to label the matrices by the indices
$0,1,2,\ldots n$ and in the solvable case by $1,2,\ldots, n$.  Then
$R^{(s)}$ will be obtained from $W^{(s)}$.

On the other hand, if we have the $n\times n$ matrices $W^{(s)}$,
$s=0,1,\ldots n-1$ we can introduce the $(n+1)\times (n+1)$ matrices
having the block form:
\begin{equation}
\begin{array}{l}
Q^{(0)}={\bf 1}_{n+1}\\[4pt]
Q^{(s)}=\left(\begin{array}{cc}0&W^{(s-1)}\\e_s&0  \end{array}\right),
\quad 1\leq s\leq n
\end{array}
\end{equation}
where $\{e_i\}_1^n$ is the canonical basis in ${\bf K}^n$. In other words,
to $W^{(s)}$, $s\geq 1$ we append zero first row and first column
equal to $e_s$. Adding to the set the unit matrix, we get a tensor
$Q_{ij}^k$, $0\leq i,j,k\leq n$ (by definition $(R^{(i)})_k^j\equiv
Q^{ij}_k$) and this tensor corresponds to a semisimple extension.
Thus the semisimple and solvable extensions are in one-to-one
correspondence and it is sufficient to study only one of these classes.

Below we shall assume that the matrices $W^{(i)}$ are already put into the
canonical form and we have the solvable case. The vector space
${\cal G}^n$ can be split into
\begin{equation}
{\cal G}^n={\cal F}_n^{(1)}\oplus{\cal F}_n^{(2)}\oplus\ldots
\oplus {\cal F}_n^{(n)}
\end{equation}
where ${\cal F}_n^{(i)}$ consists of those ${\bf x}$ for which
$x_j=0$ for all $j\neq i$. For $1\leq k\leq n$ let us set 
\begin{equation}
{\cal F}[n,k]={\cal F}_n^{(k)}\oplus {\cal F}_n^{(k+1)}\oplus\ldots 
\oplus {\cal F}_n^{(n)}
\end{equation}
and for $k>n$ let ${\cal F}[n,k]=0$. Obviously,
\begin{equation}
0\subset {\cal F}^{(n)}_n={\cal F}[n,n]\subset {\cal F}[n,n-1]\subset
\ldots\subset {\cal F}[n,1]={\cal G}_W^n
\end{equation}
and
\begin{equation}
\left[{\cal F}[n,k],{\cal F}[n,s]\right]\subset {\cal F}[n,(\max(k,s)) +1]
\end{equation}
In particular ${\cal F}_n^{(n)}$ is Abelian ideal. Denote by ${\cal
S}^{k}_{n}$ the maximal Abelian ideal of the type ${\cal
F}[n,n-k+1]$. As ${\cal F}[n,n]$ is an Abelian ideal ${\cal
S}^{n}_{k}$ is not $0$. Then we have the exact sequence
\begin{equation}
0\mapsto {\cal S}_{n}^{k}\mapsto {\cal G}^n_W\mapsto {\cal G}^n_W\slash
{\cal S}^{n}_{k}\mapsto 0
\end{equation}
which defines an extension through an Abelian kernel, see \cite{GoGr,ChEi},
in the sense that it is usually understood.  In this manner the
things are reduced to the study of ${\cal G}^{n-k}$ and the tensor
$\bar{W}$ on it, a process which we call a reduction. It amounts simply to take
the tensor $W_{k}^{ij}$ and let the indices to run over $1,2,\ldots
n-k$.  As the above reduction can be performed by $k$ reductions of
the type ${\cal G}_W^{n}\mapsto {\cal G}_{\bar{W}}^{n-1}$ we shall
always assume that $k=1$.

Now consider the case ${\cal G}_W^{n+1}$. In order to simplify the
notations we also write ${\cal G}_a$ instead of ${\cal F}[n+1,n+1]$ and
${\cal G}^{n}_{\bar {W}}$ instead of ${\cal G}^{n+1}_W\slash {\cal
F}[n+1,n+1]$.  Evidently ${\cal G}_a$ as vector space is equivalent to
${\cal G}$ but with a structure of an Abelian algebra.  It is not
difficult to see that the representation induced over ${\cal G}_a$
for the extensions of the type:
\begin{equation}\label{eq:EXT2}
0\mapsto {\cal G}_a\mapsto {\cal G}^{n+1}_W\mapsto {\cal
G}^{(n)}_{\bar{W}}\mapsto 0
\end{equation}
is the trivial one. Now reminding that to every extension with Abelian
kernel there correspond 2-cohomologies,  \cite{GoGr} it is
interesting to ask what is the space $H^2({\cal
G}^{n}_{\bar{W}},{\cal G}_a)$. (Here the notations are standard
ones, the cohomologies correspond to the algebra ${\cal
G}^{n}_{\bar{W}}$ acting trivially on ${\cal G}_a$). It is almost
readily seen,\cite{JLT} that as a vector space $H^2({\cal
G}^{n}_{\bar{W}},{\cal G}_a)$ is spanned  by the last rows of the
matrices $W^{(i)}$ and the coboundaries are the linear combinations
of the first $n$ rows of these matrices. For this reason it is
convenient to introduce the matrices $W_{(k)}$:
\begin{equation}
(W_{(k)})^{ij}=W_k^{ij}
\end{equation}
Then the corresponding element from the cohomology space
$H^2({\cal G}^{n}_{\bar{W}},{\cal G}_a)$ can be identified with
$W_{(n+1)}$ and the coboundaries are linear combinations of the
matrices $\{W_{(i)}\}_{i=1}^{n}$. Note that the $(n+1)\times (n+1)$ matrix
$W_{(n+1)}$ is symmetric and the elements in its last row and column
are equal to zero. In the same manner the elements in the rows with
numbers $j\geq i$ of $W_{(i)}$ are zero. For the above reasons,
dropping out the zero rows and columns we can
identify $W_{(i)}$ with $(i-1)\times (i-1)$ symmetric matrix.
Therefore it follows that ${\cal G}_W^n$ is also
uniquely defined by the set of symmetric matrices:
\begin{equation}
(W_{(1)},W_{(2)},\ldots, W_{(n)}),
\end{equation}
the matrix $W_{(i)}$ ($i>1)$ being of the type $(i-1)\times (i-1)$ and always
$W_{(1)}=0$. In these notations is written the table of all solvable
extensions up to $n\leq 4$ in \cite{JLT}.

\begin{Example}
For arbitrary $\lambda=(\lambda_1,\lambda_2)\in {\bf K}^2$ 
the formulae:
\begin{equation}
\begin{array}{l}
[{\bf x},{\bf y}]^{(\lambda)}_0=[x_0,y_0]\\[4pt]
[{\bf x},{\bf y}]^{(\lambda)}_1=[x_0,y_1]+[x_1,y_0]\\[4pt]
[{\bf x},{\bf y}]^{(\lambda)}_2=[x_0,y_2]+[x_2,y_0]+
\lambda_1([x_1,y_1])\\[4pt]
[{\bf x},{\bf y}]^{(\lambda)}_3=[x_0,y_3]+[x_3,y_0]+
\lambda_1([x_1,y_2]+[x_2,y_1])\\[4pt]
\dotfill\\[4pt]
[{\bf x},{\bf y}]_{n-1}^{(\lambda)}=[x_0,y_{n-1}]+
[x_{n-1},y_0]+\lambda_1([x_1,y_{n-2}]+\ldots [x_{n-2},y_1])\\[4pt]
[{\bf x},{\bf y}]^{(\lambda)}_{n}=[x_0,y_{n}]+[x_{n},y_0]+
\lambda_2([x_1,y_{n-1}]+\ldots [x_{n-1},y_1])
\end{array}
\end{equation}
define Lie algebra structure 
$({\bf x},{\bf y})\mapsto [{\bf x},{\bf y}]^{(\lambda)}$ on
${\cal G}^{n+1}$

The interpretation of this fact is the following. Considering the
semisimple part of this extension, we see that it corresponds to the
so called Leibnitz extension bracket, defined by the tensor:
\begin{equation}
\hat{W}^{ij}_k=\lambda_1\delta^{i+j}_k,\quad 1\leq i,j,k\leq n-1
\end{equation}
extended using the cocycle $R^{(n)}$: $(R^{(n)})_i^j=\lambda_2\delta_n^{i+j}$
where $\delta_m^s$ is the Kronecker delta.

The Lie algebra structure, responsible for the Poisson structure of
the so called Compressible Reduced Magnetohydrodynamics model,
\cite{JLT}is a particular case of this construction. Indeed, in this
case the nonzero components of the $W_i^{jk}$ $i,j,k=0,1,2,3$
are:
\begin{equation}
\begin{array}{c}
W_{i}^{0i}=W_{i}^{0i}=1,\quad i=0,1,2,3\\[4pt]
W_3^{12}=W_3^{21}=-\beta_e
\end{array}
\end{equation}
where $\beta_e$ (the electron $\beta$) is a parameter. The bracket
corresponding $W_i^{jk}$ is given by
\begin{equation}
\begin{array}{l}
[{\bf x},{\bf y}]_0=[x_0,y_0]\\[4pt]
[{\bf x},{\bf y}]_1=[x_0,y_1]+[x_1,y_0]\\[4pt]
[{\bf x},{\bf y}]_2=[x_0,y_2]+[x_2,y_0]\\[4pt]
[{\bf x},{\bf y}]_3=[x_0,y_3]+[x_3,y_0]-\beta_e [x_1,y_2]-\beta_e[x_2,y_1]
\end{array}
\end{equation}
and it is exactly the structure we considered is the case $n=3$;
$\lambda_1=0,\lambda_2 =-\beta_e$.
\end{Example}

\section{The universal extensions via commutative\\ algebras structures}

All the results introduced in the above are obtained via the theory of
Lie algebras and the Lie algebra extensions. However, another approach
to the universal Lie algebra extensions is also possible. As a matter of
fact, a more scrutinized look permits to perceive that an additional
algebraic structure is suggested by the properties of the
tensors $W^{ij}_k$.  Let ${\cal A}^n$ be a vector space with a basis
$\{e^i\}_{i=1}^n$ over the field ${\bf K}$.  Let us put:
\begin{equation}\label{eq:glt}
e^i*e^j=\sum\limits_{s=1}^n W^{ij}_s e^s
\end{equation}
where $W^{ij}_k$ is a tensor with the properties (\ref{eq:0}).
Let us extend this product by bi-linearity over ${\cal A}^n$. Then
it is not difficult to see that ${\cal A}^n$ becomes an associative
commutative algebra which we shall denote by ${\cal A}_W^n$.
Conversely, if ${\cal B}^n$ is an associative commutative
finite-dimensional algebra, then by (\ref{eq:glt}) we can define the
set $W^{ij}_k$ which will possess the properties (\ref{eq:0}).  Thus
we obtain:
\begin{Theorem}\label{Theorem:1}
The $n$-dimensional universal extensions ${\cal G}^n_W$ over ${\bf
K}$ defined by the tensors $W^{ij}_k$ are in one-to-one
correspondence with the $n$-dimensional associative commutative
algebras ${\cal A}_W^n$ over the same field. Moreover, the correspondence
\begin{equation}
e^i\mapsto W^{(i)}\in \Mat(n, {\bf K}),\quad (W^{(i)})^j_k=W^{ij}_k,
\end{equation}
defines a matrix representation of the algebra ${\cal A}_W^n$ and hence
to any linear transformation of the elements $e^i$
\begin{equation}
e^{i'}=\sum\limits_{i=1}^n A^{i'}_i e^i
\end{equation}
correspond the transformations (\ref{eq:trm1}) and (\ref{eq:trm2})
of the tensor $W^{ij}_k$ and of the matrices $W^{(i)}$.
\end{Theorem}
Only the second part of this theorem is not obvious, but the proof is
obtained by a simple computation. We want only to note that $\Mat(n,
{\bf K})$ is considered here as an associative algebra and that by
representation we understand a linear map $\Phi: {\cal A}^n\mapsto
\Mat(m, {\bf K})$, such that:
\begin{equation}
\Phi(x*y)=\Phi(x)\Phi(y),\quad x,y\in {\cal A}^n.
\end{equation}
From the above result it is patent that the matrix $W^{(i)}$ is in
fact the matrix of the action of the element $e^i$ on ${\cal A}_W^n$.
This permits to understand easily the block structure and the
canonical structure of the matrices $W^{(i)}$:
\begin{Proposition}
The splitting of ${\cal A}_W^n$ into a sum of ideals:
\begin{equation}
{\cal A}_W^n={\cal J}^{s_1}_{W_1}\oplus {\cal J}^{s_2}_{W_2}\oplus
\ldots\oplus{\cal J}^{s_q}_{W_q} 
\end{equation}
with dimensions $s_i=\dim{\cal J}^{s_i}_{W_i}$ corresponds to the
splitting of the algebra ${\cal G}_W$ into a direct sum:
\begin{equation}
{\cal G}_W^n={\cal G}^{s_1}_{W_1}\oplus{\cal
G}^{s_2}_{W_2}\oplus\ldots \oplus{\cal G}^{s_q}_{W_q}
\end{equation}
The irreducible ${\cal G}^{n}_{W}$ corresponds to ${\cal A}_W^n$ that
cannot be split into a sum of proper ideals, that is to simple
${\cal A}_W^n$. 
\end{Proposition}
\begin{Proposition}
The canonical structure of $W^{(i)}$ is equivalent to the requirement
that for the solvable case in the basis $\{e^i\}_{i=1}^n$  
\begin{equation}
\begin{array}{lr}
e^i*e^j=
\sum\limits_{s>max(i,j)}W^{ij}_s e^s&1\leq i,j\leq n
\end{array}
\end{equation}
and in the semisimple case we have one more independent element $e^0$
such that 
\begin{equation}
\begin{array}{lr}
e^0*e^j=e^j,&0\leq j\leq n\\
\end{array}
\end{equation}
\end{Proposition}
Thus the solvable case corresponds to the situation when
the operators $T_a$, $a\in {\cal A}^n_W$
\begin{equation}
T_a(b)\equiv a*b
\end{equation}
are nilpotent and the semisimple case corresponds to algebras with unity.

We shall continue to consider only the solvable case ${\cal A}_W^n$.
In order to pass from the semisimple case to the solvable one simply
must put away the unity element and in order to pass from the
solvable to the semisimple one must add an unity element.

Let ${\cal A}_n^{(s)}$ be the
one-dimensional vector space spanned by $e^s$ and for $1\leq k\leq n$
let 
\begin{equation}
{\cal A}[n,k]={\cal A}_n^{(k)}\oplus {\cal A}_n^{(k+1)}\oplus\ldots 
\oplus {\cal A}_n^{(n)}
\end{equation}
For $k>n$ we put again ${\cal A}[n,k]=0$. Then
\begin{equation}
0\subset {\cal A}^{(n)}_n={\cal A}[n,n]\subset {\cal A}[n,n-1]\subset
\ldots\subset {\cal A}[n,1]={\cal A}_W^n
\end{equation}
Clearly,
\begin{equation}
{\cal A}[n,k]*{\cal A}[n,s]\subset {\cal A}[n,(\max(k,s))+1]
\end{equation}
for arbitrary $1\leq k,s\leq n$.

Let us turn now our attention to the extensions from ${\cal A}_W^n$ to
an algebra ${\cal A}^{n+1}$ and assume that the matrices $W^{(i)}$
can be put into canonical form. Then the construction of the extension
we are looking for consists in adding to the basis $\{e^i\}_{i=1}^n$
the vector $e^{n+1}$ in such a way that ${\cal J}={\bf K}e^{n+1}$ is
an ideal and the sequence
\begin{equation}
0\mapsto {\cal J}\mapsto {\cal A}^{n+1}\mapsto {\cal A}_W^{n}\mapsto 0
\end{equation}
is exact.
Assuming that the basis $\{e_i\}_{i=1}^{n+1}$ in ${\cal A}^{n+1}$ is
also canonical and denoting the multiplication in it by a
point, we have:
\begin{equation}
\begin{array}{lr}
e^i.e^j=\sum\limits_{k=1}^nW^{ij}_ke^k+R^{ij}e^{n+1},&
1\leq i,j\leq n\\[8pt]
e^i.e^{n+1}=0,&1\leq i\leq n+1\\[8pt]

\end{array}
\end{equation}
$R^{ij}=R^{ji}$ and it clearly plays the role of $W_{n+1}^{ij}$ in the
algebra ${\cal A}^{n+1}$. Next, it is
easily seen that in order to ensure the associative rule, $R^{ij}$
must satisfy:
\begin{equation}\label{eq:01}
\sum\limits_{s=1}^n(R^{is}W_{s}^{jk}-R^{js}W_{s}^{ik})=0
\end{equation}
Thus $R^{ij}$ satisfying (\ref{eq:01}) determines the extension. We
remind that $R^{ij}$ are just the cohomologies defined by the
universal Lie algebra extensions, and 
\begin{equation}
\bar{R}^{ij}=\sum\limits_{s=1}^n\lambda^s W^{ij}_s,\quad \lambda^s\in
{\bf K}
\end{equation}
are the corresponding coboundaries.

From the other side, we can say that every element $\xi$ in ${\cal
A}^{n+1}$ can be written uniquely into the form
\begin{equation}
\xi={\bf x}+xe^{n+1}=\sum\limits_{k=1}^n x_ie^i+xe^{n+1}
\end{equation}
and if $\eta={\bf y}+ye^{n+1}$ then
\begin{equation}
\xi.\eta={\bf x}*{\bf y}+R({\bf x},{\bf y})e^{n+1}
\end{equation}
where
\begin{equation}
R({\bf x},{\bf y})=\sum\limits_{k,s=1}^n R^{ij}x_iy_j
\end{equation}
is a symmetric bilinear function on ${\cal A}_W^n$ with values in ${\bf
K}$. The condition (\ref{eq:01}) means that
\begin{equation}
R({\bf x}*{\bf y},{\bf z})=R({\bf x},{\bf y}*{\bf z})
\end{equation}
If we introduce in $({\cal A}^n)^*$ the dual basis $\{e_s\}_{s=1}^n$ 
($\langle e^j,e_i\rangle=\delta_i^j$) then the coboundaries can be regarded as
bilinear functions of the type:
\begin{equation}
r({\bf x},{\bf y})=\sum\limits_{s=1}^n\lambda^s\langle{\bf
x}*{\bf y},e_s\rangle
\end{equation}

We are able now to introduce cohomologies related to ${\cal A}^{n}$ in
which the above objects, namely $R$ and $r$, have the same
meaning as they have with respect to the Lie algebra structures. To this end
let us introduce the following complex:
\begin{itemize}
\item The cochains. For $s=0$ we set
$C^0({\cal A}^n_W)={\bf K}$ and for $s\geq 1$, let $C^s({\cal A}_W^n)$ be
the set of all the polylinear functions $\omega_s$
\begin{equation}
\omega_s:{\cal A}^n\mapsto {\bf K}.
\end{equation}
\item The coboundary operator $d_p$. For $\lambda\in C^0({\cal
A}_W^n)={\bf K}$ we put $d_0(\lambda)(a)=0$, and for $\omega_p\in
C^p({\cal A}_W^n), \quad p\geq 1$ we define $d_p\omega_p\in
C^{p+1}({\cal A}_W^n)$ by the formula:
\begin{equation}\label{eq:ch}
\begin{array}{l}
d_p\omega_p(a_1,a_2,\ldots,a_{p+1})=\\[4pt]
\sum\limits_{k=1}^p(-1)^k
\omega_n(a_1,a_2,\ldots,a_k*a_{k+1},a_{k+2},\ldots,a_{p+1})
\end{array}
\end{equation}
$a_i\in {\cal A}^n,\quad 1\leq i\leq p+1$.
\end{itemize}
One can check that $d_{p}\circ d_{p+1}=0$ and therefore the above
construction defines a complex. The resulting cohomology groups we
shall denote by $H^p({\cal A}_W^n)$. Similar notations we
use also for the cocycles and coboundaries, that is they will be
denoted by $Z^p({\cal A}_W^n)$ and $B^p({\cal A}_W^n)$.

Now it is clear that $R({\bf x},{\bf y})$ defines a cohomology class in
$H^2({\cal A}_W^n)$ and $r({\bf x},{\bf y})$ is again a
coboundary. Indeed, if $\alpha: {\cal A}^n\mapsto {\bf K}$ is a
linear map, then 
$$
[d_1\alpha]({\bf x},{\bf y})=-\langle{\bf x}*{\bf y},\alpha\rangle
$$ 
whence
\begin{equation}
R=d_1\beta,\quad \beta=-\sum\limits_{s=1}^n\lambda^s e_s
\end{equation}
Thus even as regards to the cohomologies the correspondence between
the Lie algebras ${\cal G}_W^n$ and the commutative algebras ${\cal
A}^{n}$ is complete, that is the spaces $H^2({\cal A}^n_W)$ and
$H^2({\cal G}_W,{\cal G}_a)$ coincide.

The structure of a vector space on ${\cal A}^n$ seems indispensable,
but in fact there is a special case when we can manage with weaker
structure. In order to show this, let ${\cal
M}=\{e_1,e_2,\ldots,e_n\}$ be only a commutative monoid with
multiplication operation $*$. Let us denote
$\bar{I}_n=\{1,2,...,n\}$. Then
\begin{equation}
e^i*e^j=e^{f(i,j)},
\end{equation}
where $f$ is a function from $\bar{I}_n\times \bar{I}_n$ to
$\bar{I}_n$, such that $f(i,f(j,k))=f(f(i,j),k)$ and $f(i,j)=f(j,i)$.
We shall call such a function an $E$-function.

It is readily seen that it is enough to put
$W^{ij}_k=\delta^{f(i,j)}_k$ in order to obtain a universal extension tensor.
In this way we actually construct a finite dimensional commutative subalgebra
of $\Mat(n,{\bf K})$ for arbitrary ${\bf K}$ because the elements of
$W^{(i)}$ are either $0$ or $1$. Then we can impose the vector space
structure, assuming that $\{e^i\}_{i=1}^n$ are independent ''generators'' and
come to the corresponding algebra ${\cal A}^n$ which is of course the
free algebra corresponding to ${\cal M}$.

Let us see how the outlined construction works in the case of two
very simple monoid structures, the ones defined by the addition and
the multiplication in the ring ${\bf Z}_p$ - the ring of integers
$\mod (p)$.  Taking the addition structure we easily get:
\begin{equation}
W^{ij}_k=\delta_k^{i+j},\quad i,j=0,1,2\ldots, p-1,
\end{equation}
where the sum is understood modulo $p$. This turns out to be a
splitting extension, see \cite{ABY}, that is ${\cal G}^p_W={\cal G}^p$. It
can be shown however, that with a suitable deformation, see \cite{ABY}, it can
be transformed into the so called Leibnitz extension. (For its
definition see \cite{JLT} or the next section).

Let us take now the multiplicative structure. The case $p=2$ being
trivial, let us take $p=3$. Then the matrices $W^{(i)}$ are:
\begin{equation}
\begin{array}{ccc}
W^{(0)}=\left(\begin{array}{ccc}1&1&1\\0&0&0\\0&0&0 \end{array}\right) &
W^{(1)}=\left(\begin{array}{ccc}1&0&0\\0&1&0\\0&0&1 \end{array}\right) &
W^{(2)}=\left(\begin{array}{ccc}1&0&0\\0&0&1\\0&1&0 \end{array}\right)
\end{array}
\end{equation}
If we regard $W^{(i)}$ as corresponding to the independent generators
$e^i$, then it is easily seen that the linear transform:
\begin{equation}
\begin{array}{l}
n^0=e^0\\
n^1=-e^0+{1\over 2}(e^1+e^2)\\
n^2={1\over 2}(e^1-e^2)
\end{array}
\end{equation}
changes $W^{(i)}$, $i=0,1,2$, into the set
\begin{equation}
\begin{array}{ccc}
N^{(0)}=\left(\begin{array}{ccc}1&0&0\\0&0&0\\0&0&0 \end{array}\right) &
N^{(1)}=\left(\begin{array}{ccc}0&0&0\\0&1&0\\0&0&0 \end{array}\right) &
N^{(2)}=\left(\begin{array}{ccc}0&0&0\\0&0&0\\0&0&1 \end{array}\right)
\end{array}
\end{equation}
and hence the corresponding extension splits: 
${\cal G}^3_W={\cal G}\oplus {\cal G}\oplus {\cal G}$.

Let us consider $p=4$. Instead of
writing down the matrices, let us start directly by the the
generators $e^i$, $i=0,1,2,3$, $e^i*e^j=e^{ij(mod(4))}$. 
If we introduce $n^0=e^0$, $n^i=e^i-e^0$, then
$n^i*n^0=\delta^{i0} n^i$ and $n^i*n^j=n^{ij(mod(4))}$ for
$i,j=1,2,3$.  The matrices in this basis have block-diagonal
form, that is the extension splits : ${\cal G}^4_W={\cal
G}\oplus{\cal G}^3_N$ where $N_{k}^{ij}$, $1\leq i,j,k\leq 3$ is a
new tensor describing the structure in ${\cal G}^3_N$. It is defined by
the matrices $N^{(i)}$:
\begin{equation}
\begin{array}{ccc}
N^{(1)}=\left(\begin{array}{ccc}1&0&0\\0&1&0\\0&0&1 \end{array}\right) &
N^{(2)}=\left(\begin{array}{ccc}0&0&0\\1&0&1\\0&0&0 
\end{array}\right)&
N^{(3)}=\left(\begin{array}{ccc}0&0&1\\0&1&0\\1&0&0 
\end{array}\right)
\end{array}
\end{equation}
The change:
\begin{equation}
\begin{array}{l}
f^1=e^1,\quad f^2=e^1-e^3,\quad f^3=e^2
\end{array}
\end{equation}
transforms the matrices $N^{(i)}$ into the set of the matrices $F^{(i)}$:
\begin{equation}
\begin{array}{ccc}
F^{(1)}=\left(\begin{array}{ccc}1&0&0\\0&1&0\\0&0&1 \end{array}\right) &
F^{(2)}=\left(\begin{array}{ccc}0&0&0\\1&1&0\\0&0&0 
\end{array}\right)&
F^{(3)}=\left(\begin{array}{ccc}0&0&0\\0&0&0\\1&0&0 
\end{array}\right)
\end{array}
\end{equation}
Then putting away the semisimple part we get:
\begin{equation}
\begin{array}{cc}
M^{(1)}=\left(\begin{array}{cc}1&0\\0&0 \end{array}\right) &
M^{(2)}=\left(\begin{array}{cc}0&0\\0&0 
\end{array}\right)
\end{array}
\end{equation}
which corresponds to a extension with $n=2$ that splits.
More interesting examples of the described type will be exhibited in the 
next section. 

\section{The commutative structures corresponding to the irreducible
solvable extensions}

In the Introduction mentioned that in \cite{JLT} are classified the
irreducible solvable extensions up to $n=4$ over ${\bf C}$ and is
presented a table with the matrices $W_{(i})$ from which one can
construct the extension. If we peruse all the cases in the table we can remark
the following interesting fact: for every fixed case with tensor
$W^{ij}_k$, $1\leq i,j,k\leq n$, there exists a function:
\begin{equation}
\begin{array}{c}
f:I_n\times I_n\mapsto I_n\\[4pt]
I_n=\{0,1,2,\ldots,n\}
\end{array}
\end{equation}
such, that 
\begin{equation}\label{eq:CF2}
W^{ij}_k=\delta_k^{f(i,j)},\quad 1\leq i,j,k\leq n
\end{equation}
and possessing the properties:
\begin{enumerate}
\item[I.] $f$ is symmetric
\begin{equation}\label{eq:an1}
f(i,j)=f(j,i)
\end{equation}
\item[II.] $f$ satisfies:
\begin{equation}\label{eq:an2}
f(i,f(j,k))=f(f(i,j),k), \quad i,j,k\in I_n
\end{equation}
\item[III.] For all $i$:
\begin{equation}\label{eq:an3}
f(i,0)=0
\end{equation}
\item[IV.] If $f(i,j)\neq 0$ then $f(i,j)>\max(i,j)$. In particular, this 
means that $f(i,n)=0$ for every $n$.
\end{enumerate}
Let us remark, that if as before $\bar{I}_n=\{1,2,\ldots,n\}$, 
then the function $f$ is equal to zero on $D_n\equiv (I_n\times
I_n)\setminus (\bar{I}_n\times \bar{I}_n)$, that is on the
''boundary'' of the square $(I_n\times I_n)$.

We see that at least for $n\leq 4$ to the irreducible nonequivalent
solvable extensions correspond specific monoids that can be described
in the following way:
{\it Suppose that the set ${\cal M}$ contains $n+1$ different elements 
$${\cal M}=\{Z^0,Z^1,{Z}^2,\ldots, Z^n\}$$ and let us introduce in
${\cal M}$ the operation :
\begin{equation}\label{eq:BO}
Z^i*Z^j=Z^{f(i,j)}
\end{equation}
where $f(i,j)$ satisfies the properties I-IV.
Then the set ${\cal M}^0$ with
the binary operation (\ref{eq:BO}) defines a commutative monoid such
that:}
\begin{equation}
\begin{array}{c}
Z^0*Z^i=Z^0 \quad \mbox{ for all i}\\
\mbox{if} \quad Z^i*Z^j=Z^k\neq Z^0\quad \mbox{then}\quad
 k>\max(i,j).
\end{array}
\end{equation}
The relation of this structure with the structure described in
Theorem (\ref{Theorem:1}) is the following.  Suppose that
we already introduced the algebra ${\cal A}^{n+1}$ corresponding to
${\cal M}^0$ by the $n+1$ basis vectors $\{e^i\}_{0}^{n}$,
such that $e^i*e^j=e^{f(i,j)}$. Then ${\cal J}={\bf K}e^0$ is an ideal and
$\bar{\cal A}^{n}\equiv {\cal A}^{n+1}\slash {\cal J}$ is also an
algebra. As $e^i$, $i\neq 0$ can be taken as representatives in the
equivalence classes the multiplication in $\bar{\cal A}^{n}$ it is given by:
\begin{equation}
e^i*e^j=\left\{\begin{array}{lr} e^{f(i,j)}&f(i,j)\neq 0\\
0&f(i,j)=0
\end{array}\right.
\end{equation}
Therefore the algebra that corresponds to the solvable extension under
consideration is $\bar{\cal A}^n$.

Clearly, (\ref{eq:CF2}) with $f$ satisfying the requirements I-IV
always defines a solvable extension, but whether a general 
solvable $W_k^{ij}$ with a linear transformation can be put into this
form is an open question.

The functions $f(i,j)$ with the properties I-IV seem to
play an important role in the theory. We shall call them
$SE$-functions. As already mentioned, the function $f(i,j)$ satisfies
only I-II we call $E$-function.

Let us consider some special cases of $SE$ functions. To this end let
$I_n$ and $\bar{I}_n$ be as above, that is $I_n=\{0,1,2,\ldots,n\}$ and
$\bar{I}_n=\{1,2,\ldots,n\}$.
\begin{itemize}
\item 
The Leibnitz extension, see \cite{JLT} of dimension $n$ is characterized by
\begin{equation}
f_L(i,j)=\left\{\begin{array}{lcc}i+j&\mbox{if}& i+j<n \\0&\mbox{if}&i+j\geq n
\end{array}\right.
\end{equation}
when $1\leq i,j\leq n$ and $f_L(0,i)=0$.
\item Extension from ${\cal G}_a^{k}$-Abelian to ${\cal
G}_f^n$, $n>k$ :
\begin{equation}
f(i,j)=\left\{\begin{array}{lcl} n&\mbox{if}& (i,j)\in A\\
 0&\mbox{if}& (i,j)\in (I_n\times I_n)\setminus A 
\end{array}\right.
\end{equation}
Here $A$ is arbitrary symmetric subset of $\bar{I}_n\times\bar{I}_n$.
(We  call a set $A\subset I_n\times I_n$ symmetric if from $(s,l)\in
A$ follows that $(s,l)\in A$.)

\item Trivial extension from ${\cal G}^{k}$, characterized by function
$\hat{f}$, to ${\cal G}^n$, ($k<n$):
\begin{equation}
f(i,j)=\left\{\begin{array}{lcl}\hat{f}(i,j)&\mbox{if}&(i,j)\in
{I}_k\times I_k 
\\0&\mbox{if}&(i,j)\in (I_n\times I_n)\setminus (I_k\times I_k)
\end{array}\right.
\end{equation}
\end{itemize}
More generally, let $f$ be an $SE$-function defined in $I_n\times I_n$ and let
\begin{equation}
A_f=f^{-1}(n)=\{(i,j):f(i,j)=n\}\subset I_{n-1}\times I_{n-1}.
\end{equation}
For $(i,j)\in I_{n-1}\times I_{n-1}$ we set:
\begin{equation}
\hat{f}(i,j)=\left\{ \begin{array}{lcl}
f(i,j)&\mbox{if}& (i,j)\in
({I}_{n-1}\times {I}_{n-1})\setminus A_f\\
0&\mbox{if}& (i,j)\in A_f
\end{array}
\right.
\end{equation}
The function $\hat{f}$ is again $SE$-function and we shall say that
$f$ is obtained by extension from $\hat{f}$ or that $\hat{f}$ is a
restriction of $f$.
The ''restriction'' from $f$ to $\hat{f}$ is simply described 
through the monoid operation. Indeed, it is clear that $A_f$ is the
set of pares $(i,j)$ such, that $Z^i*Z^j=Z^n$. If we define a new product
$*_r$ by the rules:
\begin{equation}
\begin{array}{lcl}
Z^i*_r Z^j=Z^i*Z^j&\mbox{if}& 0\leq i,j\leq n-1,~~\mbox{and}~~ (i,j)\notin A_f
\\
Z^i*_r Z^j=Z^0&\mbox{if}& (i,j)\in A_f
\end{array}
\end{equation}
then it is a simple matter to see that we obtain a new monoid 
generated by the elements $Z^0,Z^1,\ldots, Z^{n-1}$.

\section{Conclusion}
We have considered the universal extensions and have shown that they
are in fact equivalent to finite dimensional commutative algebras as
all the constructions involving the universal extensions can be
reformulated as constructions for the corresponding commutative
structures. We believe that the above fact can help to understand
better both the Lie algebra and commutative algebra structures and to
stimulate further developments with a direct impact to the study of
the Poisson-Lie structures of different physical theories.


\newpage
\begin{thebibliography}{R99} 
\bibitem{JLT} J. L. Thiffeault, P. J. Morrison, Classification and
Casimir invariants of Lie-Poisson brackets, Physica D, {\bf 136}, (3-4),
p. 205-244 (2000); Jean-Luc Thiffeault, PhD Thesis, The University of Texas
at Austin, USA, December, 1998, math-ph/0009017, 11 Sep 2000
\bibitem{JLT1} J. L. Thiffeault, P. J. Morrison, The twisted top,
Phys.Lett. A, {\bf 283}, (5-6), p. 335-341 (2001)
\bibitem{Ya3} A. B. Yanovski, Linear Bundles of Lie brackets and their
Applications, J. Math. Phys., {\bf 41}, N 11, p. 7869-7882 (2000),
\bibitem{TrFo} V. V. Trofimov, A. T. Fomenko, Algebra and Geometry of 
the Integrable Hamiltonian Differential Equations, Moscow, Fizmatlit
Publishing Company, (1992), 
Minsk
\bibitem{Ki} A. A. Kirillov, Unitary representations of nilpotent Lie 
groups, Russian Math. Surveys, {\bf 17}, p. 53-104, (1962). (In Russian) 
\bibitem{Lich} A. Lichnerovich, New Geometrical Dynamics, In: Ed. A. Dold, B. 
Eckmann, Lecture notes in Mathematics, {\bf 570}, Differential geometrical 
methods in Mathematical Physics, Springer, New York, (1975) 
\bibitem{MR} J. E. Marsden, T. S. Ratiu, Introduction to Mechanics and 
Symmetry, (1994) Springer-Verlag, New York-Berlin-Heidelberg-London 
\bibitem{GoGr} M. Goto and F. Grosshans (1978) Semisimple Lie algebras, 
Lecture Notes in Pure and Applied Mathematics {\bf 38}, (New York and
Basel: M.  Dekker Inc.)
\bibitem{ChEi} C. Chevalley and S. Eilenberg, Cohomology Theory of
Lie groups and Lie algebras, Transactions of the American
Mathematical society, {\bf 63} p. 85-124 (1948)
\bibitem{ABY} A. B. Yanovski, Lie algebra extensions related with
linear bundles of Lie brackets, math.DS/0108132 (2001)

\end{thebibliography}
\end{document}